\documentclass{amsart}

\parskip=\smallskipamount

\newtheorem{theorem}{Theorem}[section]
\newtheorem{lemma}[theorem]{Lemma}

\theoremstyle{definition}
\newtheorem{definition}[theorem]{Definition}

\newtheorem{problem}[theorem]{Problem}
\newtheorem{remark}[theorem]{Remark}

\begin{document}

\title{Most real analytic Cauchy-Riemann manifolds are nonalgebraizable}

\author{Franc Forstneri\v c}
\address{Institute of Mathematics, Physics and Mechanics, 
University of Ljubljana, Jadranska 19, 1000 Ljubljana, Slovenia}
\email{franc.forstneric@fmf.uni-lj.si}
\thanks{Supported in part by Research Program  P1-0291,
Republic of Slovenia}

%
%
\date{\today} 

\begin{abstract}
We give a  simple argument to the effect that most germs of generic real analytic
Cauchy-Riemann manifolds of positive CR dimension are not holomorphically embeddable 
into a generic real algebraic CR manifold of the same real codimension in a 
finite dimensional space. In particular, most such germs are not holomorphically equivalent 
to a germ of a generic real algebraic CR manifold.
\end{abstract}

\maketitle

\section*{INTRODUCTION} 
A smooth real submanifold $M\subset \mathbf{C}^n$ in a 
complex Euclidean space is said to be a 
{\it generic Cauchy-Riemann (CR) submanifold} 
of CR dimension $m$ and codimension $d$ ($m+d=n$) 
if it is locally near every point $x\in M$ 
defined by $d$ real equations $\rho_1=0,\ldots, \rho_d=0$ 
satisfying $\partial \rho_1\wedge \ldots\wedge \partial\rho_d\ne 0$. 
(Here $\partial\rho=\sum_{j=1}^n \frac{\partial\rho}{\partial z_j}  dz_j$
and $\wedge=\wedge_\mathbf{C}$.) A germ $(M,x)$ is {\it real analytic} 
(respectively {\it real algebraic}) if it is defined
locally near $x$ by real analytic (resp.\ real algebraic) functions.
Germs $(M,x)$, $(M',x')$ are {\it holomorphically equivalent\/} 
if there exists a biholomorphic map $f\colon U\to U'$ from a neighborhood $U$ 
of $x$ onto a neighborhood $U'$ of $x'$ with $f(x)=x'$
and $f(M\cap U)=M'\cap U'$. 

Beginning with Ebenfelt [4] (1996) several authors have given examples of 
analytic CR manifolds which are not locally holomorphically equivalent 
to an algebraic one (Baouendi, Ebenfelt and Rothschild 
([1], 9.11.4), ([2], 7.2); Huang, Ji and Yau [8]). S.\ Ji studied 
the algebraization problem for real analytic strongly pseudoconvex 
hypersurfaces and established the propagation of algebraization 
for hypersurfaces with maximal Cartan-Chern-Moser rank [9], [10]. 
Gaussier and Merker studied the problem for a class of tuboids [6]. 

While it seems rather difficult to decide whether a specific analytic 
CR manifold is locally holomorphically equivalent to an algebraic one, the 
phenomenon itself is not at all surprising. 
The purpose of this note is to give a very simple argument to the effect 
that {\em most germs of generic real analytic CR manifolds 
$M\subset \mathbf{C}^n$ of positive CR dimension are not 
holomorphically embeddable into any generic real algebraic 
CR manifold $M'\subset \mathbf{C}^{n'}$ of the same codimension as $M$};
in particular, 
{\em they are not holomorphically equivalent to a germ of a real 
algebraic CR manifold in $\mathbf{C}^n$}. 
More precisely, the embeddable ones form a set of 
the first category in a suitable Baire space
(theorem 1.2). The same conclusion holds for embeddings into any 
countable union of finite dimensional families of CR manifolds 
of the same codimension. 

Our proof employs an argument from [5] 
(which essentially goes back to Poincar\'e [12])
where it was proved that most germs of real analytic 
strongly pseudoconvex hypersurfaces in $\mathbf{C}^n$ for $n>1$ 
are not holomorphically embeddable into any sphere $\sum_{j=1}^{N}|z_j|^2 =1$
(Theorem 2.2 in [5]). For embeddings into infinite dimensional 
spheres see Lempert [11].

\section{The main result}
Every germ of analytic CR manifold in $\mathbf{C}^n$ of CR dimension $m$ and 
codimension $d$, with $m+d=n$, is holomorphically equivalent to one 
of the form 
$$
	M=\{ v_j=r_j(x,y,u)\colon j=1,\ldots, d\}    \eqno(1)
$$	
where $z=x+iy\in\mathbf{C}^m$, $w=u+iv\in \mathbf{C}^d$ and 
$r=(r_1,\ldots,r_d)$ is an $\mathbf{R}^d$-valued convergent power series 
without constant and linear terms. Let $\mathcal{R}$ denote the space 
of all formal power series 
$$
	r(x,y,u)=\sum_{\alpha, \beta\in\mathbf{Z}_+^m, 
	\gamma\in \mathbf{Z}_+^d}
	c_{\alpha, \beta, \gamma}  x^\alpha y^\beta u^\gamma
	\quad (c_{\alpha,\beta,\gamma} \in \mathbf{R}^d)
$$
without constant and linear terms in $2m+d$ real variables $(x,y,u)$. 
(One could put $M$ in a Chern-Moser normal form [3],
although this will not be necessary for our purposes.)
We shall identify $r\in \mathcal{R}$ with the (formal) germ at $0\in\mathbf{C}^n$ of the 
CR manifold (1). $\mathcal{R}$ is a Fr\'echet space in the topology induced by 
the seminorms $||r||_{\alpha,\beta,\gamma}=|c_{\alpha,\beta,\gamma}|$ for all multiindices
$\alpha,\beta\in Z_{+}^{m}$, $\gamma\in Z_+^d$. The convergent power series,
representing germs of real analytic CR manifolds, form a union 
$\cup_{t>0} \mathcal{R}^t \subset \mathcal{R}$ of Banach spaces
(in fact, Banach algebras) 
$$
	\mathcal{R}^t =  \{r\in \mathcal{R}\colon ||r||_t =
	\sum |c_{\alpha,\beta,\gamma}|\cdotp t^{|\alpha|+|\beta|+|\gamma|} <+\infty\} 
								\eqno(2)
$$
with the norm $||r||_t$ ([7], p.\ 15). 
For $r\in \mathcal{R}$ and $k\in \mathbf{N}$ we denote 
by $r_k$ its the truncation (Taylor polynomial) of order $k$.
Let $\mathcal{R}_k$ be the (finite dimensional real) vector space of all 
such truncations.

\begin{definition}
A manifold $M$ (1) is {\it embeddable into an algebraic model\/}
if there exists a real algebraic CR manifold $M'\subset \mathbf{C}^{n'}$
$(n'\ge n)$ of real codimension $d$ and a holomorphic embedding 
$F=(f_1,\ldots,f_{n'}) \colon U\to \mathbf{C}^{n'}$, 
defined in an open neighborhood $U\subset \mathbf{C}^n$ of $0$, such that $F$ 
is transverse to $M'$ at $0$ and $F(M\cap U)=M'\cap F(U)$.
\end{definition}

One defines {\it formal holomorphic embeddability} of a jet (1) into 
a similar jet $M'\subset \mathbf{C}^{n'}$ defined by $v'=\rho(x',y',u')$
by requiring that the composition $\rho\circ F$ is formally
holomorphically equivalent to the jet (1) (see [5]).

\begin{theorem}
Let $t>0$ and $m\ge 1$. The set of all $r\in \mathcal{R}^t$ 
for which the germ at $0$ of the real analytic CR manifold 
$M=\{v=r(x,y,u)\}$ (1) of ${\rm CRdim}M=m$ is holomorphically embeddable into 
an algebraic model is of the first category in the Banach algebra $\mathcal{R}^t$. 
The same holds for the set of germs in $\mathcal{R}^t$ or in $\mathcal{R}$ 
which are formally holomorphically embeddable into an
algebraic model. 
\end{theorem}

\begin{proof} 
Fix the dimension $n'=m'+d \ge n=m+d$ of the target space
and denote the variables by $(z',w')$, with $z'=x'+iy'\in\mathbf{C}^{m'}$
and $w'=u'+iv'\in\mathbf{C}^d$. Every germ at $0$ of a generic 
algebraic CR manifold in $\mathbf{C}^{n'}$ of CR dimension 
$m'$ is linearly equivalent to one of the form
$$  
	A\colon \ \ \rho(z',\bar z',w',\bar w')= 
	\Im w' + \widetilde \rho(z',\bar z',w',\bar w') = 0  \eqno(3)
$$
where $\rho=(\rho_1,\ldots,\rho_d)$ is a $d$-tuple 
of real-valued polynomials and $\widetilde\rho=O(2)$ 
(i.e., it only contains terms of order $\ge 2$).
We fix a germ at $0\in\mathbf{C}^n$ of an analytic CR manifold 
$M$ of the form (1) and ask whether there exists a
germ of a holomorphic embedding 
$F=(f,g)\colon (\mathbf{C}^n,0)\to (\mathbf{C}^{n'},0)$,
with $f=(f_1,\ldots,f_{m'})$ and $g=(g_1,\ldots,g_d)$,
such that the equation
$$
	\rho(f,\bar f,g,\bar g)= 
	\Im g + \widetilde \rho(f,\bar f,g,\bar g)=0   \eqno(4)
$$
defines the germ of $M$ at $0\in\mathbf{C}^n$. (Any local holomorphic 
change of coordinates of $(\mathbf{C}^n,0)$ may be included in $F$.)
Our normalizations imply  
\begin{equation*}
\begin{split}
       T_0M=\{v=0\},  & \quad T^\mathbf{C}_0 M= T_0 M\cap i  T_0 M= \{w=0\}, \\
       T_0A=\{v'=0\}, & \quad T^\mathbf{C}_0 A = \{w'=0\}.
\end{split}
\end{equation*}
Hence $g(z,w) = Bw + \widetilde g(z,w)$ for some $B\in GL_d(\mathbf{R})$
and $\widetilde g=O(2)$. Insertion into (4) gives
$$ 
	\Im \left( Bw+\widetilde g(z,w)\right) 
	+ \widetilde \rho(f,\bar f,g,\bar g) =0.         \eqno(5)
$$
Set $g^*(z,w)=B^{-1}g(z,w)=w+B^{-1}\widetilde g(z,w)$  and
$$ 
	\widetilde\rho^*(z',\bar z',w',\bar w') = 
	B^{-1}\rho(z',\bar z',Bw',B\bar w') 
	= \Im w' +\widetilde\rho^*(z',\bar z',w',\bar w').
$$	
Multiplying (5) on the left by $B^{-1}$
we see that $M$ is also defined by 
$$
       \rho^*(f,\bar f,g^*,\bar g^*)=
       \Im g^* + \widetilde\rho^*(f,\bar f,g^*,\bar g^*)=0  
$$
where $\widetilde\rho^*=O(2)$. Thus the germ $M$ also arise 
as the preimage of the algebraic CR manifold 
$\widetilde A=\{\rho^*=0\} \subset \mathbf{C}^{n'}$ by the 
holomorphic embedding $F^*=(f,g^*)$.
This shows that it suffices to consider preimages 
of algebraic manifolds (3) by (formal) holomorphic embeddings 
$$ 
	F=(f,g),\ F(0)=0,\ \  
	g(z,w)=w+\widetilde g(z,w),\ \widetilde g=O(2).        \eqno(6)
$$
The $F$-preimage of $A$ (3) is given by  
$$
	0= \rho(f,\bar f,g,\bar g) = 
	\Im g + \widetilde\rho(f,\bar f,g,\bar g)=
	v - r'(x,y,u,v)                   \eqno(7)
$$
where $r'$ is a power series containing only terms of 
order $\ge 2$. To change the equation $v=r'(x,y,u,v)$
(7) into one of the form (1) one performs the iteration
$$ 
	v^0=0, \quad v^{j+1}=r'(x,y,u,v^{j}) \quad (j=0,1,\ldots).	
$$
In the convergent (real analytic) case this amounts 
to solving (7) on $v$ by the implicit function theorem.
The iteration  converges also on the formal level, that is, 
every coefficient of order $k$ in the power series $r$ is determined after 
at most $k$ iterations and does not change during 
subsequent iterations.

\smallskip
{\em Key observation}: If the germ at $0$ of the manifold 
$v=r(x,y,u)$ (1) is the preimage of the manifold (3) 
by a (formal) holomorphic map (6) then the coefficient of 
every monomial of order $\le k$ in 
the series for $r$ is a polynomial function of the 
coefficients (and their conjugates) of order $\le k$ in 
the series for $\rho$ and $F$, and it does not depend on the 
coefficients of order $>k$ of $\rho$ or $F$. 
\medskip

To see this it suffices to observe that all operations 
with power series without a constant term which were used in 
the process have this property (since they only involve conjugation, 
addition, multiplication, and insertion of one series into 
another, and each of these operations has the stated property). 

For a fixed $n'=m'+d$ we denote by $\mathcal{H}^{n'}$ the set of all
germs of holomorphic maps $F\colon (\mathbf{C}^n,0)\to (\mathbf{C}^{n'},0)$ 
of the form (6). For fixed $n',\nu\in\mathbf{N}$ we denote by 
$\mathcal{A}^{n',\nu}$ the set of all algebraic manifolds (3)
in $\mathbf{C}^{n'}$ defined by real polynomials $\rho=(\rho_1,\ldots,\rho_d)$
of order at most $\nu$. The corresponding spaces of truncations
are denoted by a subscript. The above observation 
amounts to the following.

\begin{lemma}
Given $k,n',\nu\in \mathbf{N}$ there exists a polynomial map 
$
	P_k \colon \mathcal{H}^{n'}_k \times \mathcal{A}_k^{n',\nu} \to \mathcal{R}_k
$
whose range contains the truncation $M_k$ of every germ 
$M$ (1) which can be formally holomorphically embedded 
in an algebraic manifold (4) of degree $\le \nu$ in $\mathbf{C}^{n'}$.
\end{lemma}

We now compare the real dimensions of the source and the target spaces.
It is easily seen that $\dim\mathcal{F}^{n'}_k \approx 2n' (k+1)^{n}$
and $\dim \mathcal{A}^{n',\nu} \approx d (\nu+1)^{n'}$, so the 
dimension of the source space is $\le C(2n'\, (k+1)^n + d (\nu+1)^{n'})$ 
for some constant $C<+\infty$ independent of $k$. On the other hand,
$\dim \mathcal{R}_k \ge ck^{2m+d}$ for some $c>0$ independent
of $k$. When $2m+d>n$ (which is the case if an only if $m>0$)
the latter dimension grows faster as $k\to+\infty$ 
(for fixed values of $n', \nu$).   
Hence for a sufficiently large $k$ the image of 
the polynomial map $P_k$ is contained in 
a union of at most countably many proper closed local real 
analytic subsets of the vector space $\mathcal{R}_k$. Since the natural
linear projection $\mathcal{R}^t\to \mathcal{R}^t_k=\mathcal{R}_k$ is surjective,
we conclude that the set of all $r\in\mathcal{R}^t$ for which 
the germ (1) can be holomorphically embedded in an algebraic
model (3) of degree $\le \nu$ in $\mathbf{C}^{n'}$ is of the 
first category in the Banach space $\mathcal{R}^t$. 
The same remains true for the countable union of these sets 
over all $n',\nu\in \mathbf{N}$. This concludes the proof of the theorem.
\end{proof}

\section{Remarks and open problems}
\begin{remark}
Clearly the proof of theorem 1.2 applies to more general families
of domains and targets; our goal here was to illustrate
a general principle without aiming at the most general results.
Similar observations in a related context have been
made recently in [6] (see especially sect.\ 8).
\end{remark}

\begin{remark}
I wish to thank P.\ Ebenfelt (private communication) 
for pointing out that, at least for hypersurfaces, the result 
can also be obtained by extending the theorem on non-embeddability 
of a generic real analytic hypersurface into a sphere 
(theorem 2.2 in [5]) to show non-embeddability into 
quadrics of any signature and then applying Webster's result 
[13] to the effect that any real algebraic hypersurface can be 
embedded into a quadric in some higher dimensional space. 
In practice this does not give a shorter proof.
\end{remark}

\begin{problem}
Let $d>1$ and $t>0$. Consider the set of all 
$r\in \mathcal{R}^t$ for which the germ at $0$ of the CR manifold 
$v=r(x,y,u)$ (1) of codimension $d$ admits a local holomorphic map 
into some algebraic strongly pseudoconvex 
{\it hypersurface} $M'\subset \mathbf{C}^N$.
Is this set of the first category in $\mathcal{R}^t$ ? 
\end{problem}

\begin{problem}
What is the answer if one replaces holomorphic embeddings 
by CR embeddings of certain smoothness class ? 
In particular, does every real analytic strongly pseudoconvex 
hypersurface in $\mathbf{C}^n$ ($n>1$) admit a local CR embedding 
of class $\mathcal{C}^1$ into an algebraic model~? Into a sphere~?
\end{problem}

\begin{problem}
Is there a propagation of holomorphic embeddability into
algebraic models, similar to [10], in a suitable class 
+of real analytic (strongly pseudoconvex) CR manifolds~?
\end{problem}

\smallskip
\it Acknowledgement. \rm
I wish to thank Peter Ebenfelt and Alexander Sukhov for 
their invaluable advice concerning the state of knowledge 
on the question considered in the paper.

\bibliographystyle{amsplain}

\end{document}